\documentclass[12pt]{amsart}

\usepackage{vmargin}
\usepackage{latexsym}
\usepackage{amsfonts}

\setpapersize{USletter}
\setmarginsrb{1.55in}{1in}{1.05in}{1in}{12pt}{32pt}{12pt}{30pt}

\begin{document}

\title{Fock space representations and crystal bases for $C_n^{(1)}$}

\thanks{{\ifcase\month\or Jan.\or Feb.\or March\or April\or May\or
June\or
July\or Aug.\or Sept.\or Oct.\or Nov.\or Dec.\fi\space \number\day,
\number\year}}

\author
	[Alejandra Premat]{{Alejandra Premat}}

\address
		{Department of Mathematics\\
		York University\\
		Toronto, Ontario\\
		Canada}
		
\email{alpremat@mathstat.yorku.ca}


\newenvironment{thm}{\subsection{}{\textbf {Theorem.}}\em}{\smallskip}
\newenvironment{prop}{\subsection{}{\textbf {Proposition.}}\em}{\smallskip}
\newenvironment{cor}{\subsection{}{\textbf {Corollary.}}\em}{\smallskip}
\newenvironment{lem}{\subsection{}{\textbf {Lemma.}}\em}{\smallskip}
\newenvironment{lemnonumber}{\noindent{\textbf {Lemma.}}\em}{\smallskip}
\newenvironment{qst}{\subsection{}{\textbf {Question.}}\em}{\smallskip}
\newenvironment{defn}{\subsection{}{\textbf {Definition.}}\em}{\smallskip}
\newenvironment{ex}{\subsection{}{\textbf {Example.}}}{\smallskip}
\newenvironment{sex}{\subsubsection{\textbf{}}{\textbf {Example.}}}{\smallskip}

\newenvironment{pf}{\noindent{\textbf {Proof.}}}{\begin{flushright}
\eop \end{flushright}\smallskip}

\newcommand\eop{{{\hfil \ensuremath \Box}}}


\newcommand\epsi{\ensuremath{\varepsilon_{i}}}
\newcommand\epsj{\ensuremath{\varepsilon_{j}}}
\newcommand\epsk{\ensuremath{\varepsilon_{k}}}
\newcommand\phii{\ensuremath{\varphi_{i}}}
\newcommand\phij{\ensuremath{\varphi_{j}}}
\newcommand\phik{\ensuremath{\varphi_{k}}}

\newcommand\etii{\ensuremath{\tilde{e}_{i}}}
\newcommand\etij{\ensuremath{\tilde{e}_{j}}}
\newcommand\etik{\ensuremath{\tilde{e}_{k}}}
\newcommand\ftii{\ensuremath{\tilde{f}_{i}}}
\newcommand\ftij{\ensuremath{\tilde{f}_{j}}}
\newcommand\ftik{\ensuremath{\tilde{f}_{k}}}
\newcommand\etiz{\ensuremath{\tilde{e}_{0}}}
\newcommand\etio{\ensuremath{\tilde{e}_{1}}}
\newcommand\etit{\ensuremath{\tilde{e}_{2}}}
\newcommand\ftiz{\ensuremath{\tilde{f}_{0}}}
\newcommand\ftio{\ensuremath{\tilde{f}_{1}}}
\newcommand\ftit{\ensuremath{\tilde{f}_{2}}}

\newcommand\my{\ensuremath{\mathbf{Y}}}
\newcommand\mx{\ensuremath{\mathbf{X}}}
\newcommand\nlao{\ensuremath{N \Lambda_{0}}}
\newcommand\mas{\ensuremath{\mbox{ and }}}

\newcommand\zii{\ensuremath{\bbZ  _\iota}}
\newcommand\ziip{\ensuremath{\bbZ  _{\iota^\prime}}}
\newcommand\ziipp{\ensuremath{\bbZ  _{\iota^{\prime \prime}}}}

\newcommand\lm{\ensuremath{\lambda}}
\newcommand\fA{\ensuremath{\mathfrak A}}
\newcommand\fB{\ensuremath{\mathfrak B}}
\newcommand\fC{\ensuremath{\mathfrak C}}
\newcommand\fD{\ensuremath{\mathfrak D}}
\newcommand\fE{\ensuremath{\mathfrak E}}
\newcommand\fF{\ensuremath{\mathfrak F}}
\newcommand\fG{\ensuremath{\mathfrak G}}
\newcommand\fH{\ensuremath{\mathfrak H}}
\newcommand\fI{\ensuremath{\mathfrak I}}
\newcommand\fJ{\ensuremath{\mathfrak J}}
\newcommand\fK{\ensuremath{\mathfrak K}}
\newcommand\fL{\ensuremath{\mathfrak L}}
\newcommand\fM{\ensuremath{\mathfrak M}}
\newcommand\fN{\ensuremath{\mathfrak N}}
\newcommand\fO{\ensuremath{\mathfrak O}}
\newcommand\fP{\ensuremath{\mathfrak P}}
\newcommand\fQ{\ensuremath{\mathfrak Q}}
\newcommand\fR{\ensuremath{\mathfrak R}}
\newcommand\fS{\ensuremath{\mathfrak S}}
\newcommand\fT{\ensuremath{\mathfrak T}}
\newcommand\fU{\ensuremath{\mathfrak U}}
\newcommand\fV{\ensuremath{\mathfrak V}}
\newcommand\fW{\ensuremath{\mathfrak W}}
\newcommand\fX{\ensuremath{\mathfrak X}}
\newcommand\fY{\ensuremath{\mathfrak Y}}
\newcommand\fZ{\ensuremath{\mathfrak Z}}

\newcommand\cA{\ensuremath{\mathcal A}}
\newcommand\cB{\ensuremath{\mathcal B}}
\newcommand\cC{\ensuremath{\mathcal C}}
\newcommand\cD{\ensuremath{\mathcal D}}
\newcommand\cE{\ensuremath{\mathcal E}}
\newcommand\cF{\ensuremath{\mathcal F}}
\newcommand\cG{\ensuremath{\mathcal G}}
\newcommand\cH{\ensuremath{\mathcal H}}
\newcommand\cI{\ensuremath{\mathcal I}}
\newcommand\cJ{\ensuremath{\mathcal J}}
\newcommand\cK{\ensuremath{\mathcal K}}
\newcommand\cL{\ensuremath{\mathcal L}}
\newcommand\cM{\ensuremath{\mathcal M}}
\newcommand\cN{\ensuremath{\mathcal N}}
\newcommand\cO{\ensuremath{\mathcal O}}
\newcommand\cP{\ensuremath{\mathcal P}}
\newcommand\cQ{\ensuremath{\mathcal Q}}
\newcommand\cR{\ensuremath{\mathcal R}}
\newcommand\cS{\ensuremath{\mathcal S}}
\newcommand\cT{\ensuremath{\mathcal T}}
\newcommand\cU{\ensuremath{\mathcal U}}
\newcommand\cV{\ensuremath{\mathcal V}}
\newcommand\cW{\ensuremath{\mathcal W}}
\newcommand\cX{\ensuremath{\mathcal X}}
\newcommand\cY{\ensuremath{\mathcal Y}}
\newcommand\cZ{\ensuremath{\mathcal Z}}

\newcommand\bbA{\ensuremath{\mathbb A}}
\newcommand\bbB{\ensuremath{\mathbb B}}
\newcommand\bbC{\ensuremath{\mathbb C}}
\newcommand\bbD{\ensuremath{\mathbb D}}
\newcommand\bbE{\ensuremath{\mathbb E}}
\newcommand\bbF{\ensuremath{\mathbb F}}
\newcommand\bbG{\ensuremath{\mathbb G}}
\newcommand\bbH{\ensuremath{\mathbb H}}
\newcommand\bbI{\ensuremath{\mathbb I}}
\newcommand\bbJ{\ensuremath{\mathbb J}}
\newcommand\bbK{\ensuremath{\mathbb K}}
\newcommand\bbL{\ensuremath{\mathbb L}}
\newcommand\bbM{\ensuremath{\mathbb M}}
\newcommand\bbN{\ensuremath{\mathbb N}}
\newcommand\bbO{\ensuremath{\mathbb O}}
\newcommand\bbP{\ensuremath{\mathbb P}}
\newcommand\bbQ{\ensuremath{\mathbb Q}}
\newcommand\bbR{\ensuremath{\mathbb R}}
\newcommand\bbS{\ensuremath{\mathbb S}}
\newcommand\bbT{\ensuremath{\mathbb T}}
\newcommand\bbU{\ensuremath{\mathbb U}}
\newcommand\bbV{\ensuremath{\mathbb V}}
\newcommand\bbW{\ensuremath{\mathbb W}}
\newcommand\bbX{\ensuremath{\mathbb X}}
\newcommand\bbY{\ensuremath{\mathbb Y}}
\newcommand\bbZ{\ensuremath{\mathbb Z}}


%
\newdimen\hoogte    \hoogte=12pt    
\newdimen\breedte   \breedte=14pt   
\newdimen\dikte     \dikte=0.5pt    
\def\beginYoung{
       \begingroup
       \def\vr{\vrule height0.8\hoogte width\dikte depth 0.2\hoogte}
       \def\fbox##1{\vbox{\offinterlineskip
                    \hrule height\dikte
                    \hbox to \breedte{\vr\hfill##1\hfill\vr}
                    \hrule height\dikte}}
       \vbox\bgroup \offinterlineskip \tabskip=-\dikte \lineskip=-\dikte
            \halign\bgroup &\fbox{##\unskip}\unskip  \crcr }
\def\End@Young{\egroup\egroup\endgroup}
\newenvironment{Young}{\beginYoung}{\End@Young}
%


\maketitle
\markboth{\textsc{Alejandra Premat}}
{\textsc{Fock space representations and crystal bases for the Affine Kac-Moody Lie algebras}}

{\allowdisplaybreaks
\abovedisplayskip 10pt
\belowdisplayskip 10pt
\jot 10 pt


\begin{center}
    \textbf{Abstract}
\end{center}
\thispagestyle{empty}
\begin{quote}
   We describe the Fock space representations and crystal bases for the affine Kac-Moody Lie algebra of type $C_n^{(1)}$ in terms of  coloured Young diagrams.  
\end{quote}

\section*{Introduction}\label{intro}

\setcounter{section}{0}
\setcounter{subsection}{0}

\bigskip
In \cite{JMMO}, the crystals for a highest weight representations of the quantized universal enveloping algebra of $\hat{\mathfrak{sl}}(n)$ are given in terms of   coloured Young diagrams  using the Fock space representations of this algebra.   In \cite{Kang}, \cite{KaKwF} and \cite{KaKwFL}, these results are generalized to some representations (those with dominant integral highest weights of level 1) of the affine Lie algebras of types $A_{2n-1}^{(2)},\ n \geq 3,\  D_{n}^{(1)},\ n \geq 4,\ A_{2n}^{(2)},\ n \geq 2,\ D_{n+1}^{(2)},\ n \geq 2,$ and $\ B_{n}^{(1)},\ n \geq 3.$  The new combinatorial objects for these algebras  are called Young walls and are built out of  cubes and ``half-cubes".   Using Young walls to obtain the results for the algebras of type  $C_n^{(1)}$ proved to be more difficult.  In \cite{KaHoCC} and \cite{KaKwFC}, the authors obtain the results for $\mathfrak{g}$ of type $C_2^{(1)}$ and in \cite{HoKaLe} they obtain a description of the crystal base for  $C_n^{(1)}$ using Young walls.    In this paper,  we define the Fock space representations of any fundamental weight for the quantized universal enveloping algebra  $\mathfrak{g}$ of type $C_n^{(1)},\ n \geq 2,$  using Young diagrams (two-dimensional combinatorial objects)  coloured appropriately - the colouring is different from that for  $\hat{\mathfrak{sl}}(n).$  We then obtain  a description  of the crystal of an irreducible representaion with a dominant highest weight in terms  of  coloured Young diagrams.  


\section*{Preliminaries}\label{Preliminaries}

\setcounter{section}{1}
\setcounter{subsection}{0}

In this section we set up the notation and state some definitions which will be needed in the following sections. 
Let $I=\{0,1,\ldots, n\}$ and
\[
A=(a_{ij})_{i,j \in I} =\left (
\begin{array}{rrrrrrr}
2 &-1&0 &0 &&\dots &\\
-2 &2&-1 &0 &&\dots &\\
0 &-1&2 &-1 &&\dots &\\
 && &\ddots &&\ &\\
 &\dots& &-1 &2&-1 &0\\
 &\dots& &0 &-1&2 &-2\\
 &\dots& &0 &0&-1 &2\\
\end{array}
\right ).
\]

  Let $(\mathfrak{h}=\bigoplus _{i=0} ^n \bbQ h_i \oplus \bbQ d, \Pi =\{h_i \ : i \in I \}, \Pi ^{\vee} =\{\alpha_{i} \ :i \in I \})$  be a realization of $A$ (see  \cite{Kac85})  and $\mathfrak{g}$ denote the affine Kac-Moody Lie algebra of type $C_n^{(1)}$ .  Hence we have $\alpha_j(h_i)=a_{ij}.$  Let  $d \in \mathfrak h$ be an element of $\mathfrak h$ such that $\alpha _i (d)=\delta_{0,i}.$  Define $s_0=s_n=2$ and $s_i=1$ for $i \in I,\ i \not = 0,n,$ so we have that $(s_i a_{ij})_{i,\ j \in I}$ is symmetric.

The \textbf{quantized universal enveloping algebra} of $\mathfrak g$, $\cU _q(\mathfrak 
g),$ is the associative algebra over $\bbQ (q)$ generated by the elements $e_i,\ f_i,\ i \in I,$ and $q^h,\ h \in P^{\vee}=\bigoplus _{i=0} ^n \bbZ h_i \oplus \bbZ d,$  subject to the following relations:

\begin{eqnarray}
&q^0 =1,\ q^h q^{h^\prime} =q^{h +h^\prime}, \mbox{ for } h \mbox{ and } h^\prime \in P^{\vee}, \\
&q^h e_i q^{-h}=q^{\alpha _i (h)} e_i, \mbox{ for } h \in P^{\vee}, \\
&q^h f_i q^{-h}=q^{-\alpha _i (h)} f_i, \mbox{ for } h \in P^{\vee}, \\
&e_if_j -f_j e_i =\delta_{ij} \frac {k_i - k_i ^{-1}}{q_i -q_i^{-1}}, \mbox{ for } i, j, \in I,\\
&\sum _{k=o}^{1-a_{ij}}\, (-1)^k  \left [ \begin{array}{c} 1-a_{ij} \\ k \end{array}\right ]_{q_i} e_i^{1-a_{ij}-k} e_j e_i^k =0 \mbox{ for } i \not = j,\\
&\sum _{k=o}^{1-a_{ij}}\, (-1)^k  \left [ \begin{array}{c} 1-a_{ij} \\ k \end{array}\right ]_{q_i} f_i^{1-a_{ij}-k} f_j f_i^k =0 \mbox{ for } i \not = j,
\end{eqnarray}

\noindent
where $q_i=q^{s_i},\ k_i =q^{s_i h_i}, \ [n]_{q_i}=\frac{q_i^n - q_i^{-n}}{q -q^{-1}},\ [0]!=1,\ [n]_{q_i} ! =[n]_{q_i}[n-1]_{q_i}\cdots[1]_{q_i}$ and $\left [ \begin{array}{c} m \\ n \end{array}\right ]_{q_i}=\frac{[m]_{q_i} ! }{[n]_{q_i}! [m-n]_{q_i} ! }.$


\noindent
The \textbf{weight lattice} for $\mathfrak{g}$ is defined to be $P:=\{ \lambda \in \mathfrak{h} ^* : \lambda (P^{\vee}) \subset \bbZ \}$.

\noindent A $\cU _q(\frak g)-$module M is said to belong to the 
category $\cO_{int}$ if
\begin{enumerate}
\item[(i)] $ \ M=\bigoplus _{\lambda \in P} M_\lambda \text { where } 
M_\lambda 
=\{u\in M : q^hu=q^{\lambda (h)}u \ \ \forall h \in 
P^\vee\}$,
\item[(ii)] $\ \mathrm{dim}(M_\lambda)<\infty \mbox{ for all } \lambda \in P $,
\item[(iii)] for each $i \in I,\ M $ is the union of finite dimensional  
$\cU _q({\frak g_i})$-modules where ${\frak g_i}$ is the subalgebra of 
$\frak g$ generated by $e_i,f_i,q^{h_{i}}\text { and }  q^{-h_{i}}$, and
\item[(iv)] $\ M=\bigoplus _{\lambda \in F+Q_-} M_\lambda$, where $F$ is a 
finite subset of $P$ and $Q_{-}=-\sum_{i \in I } \bbN \alpha _i.$
\end{enumerate}
\noindent The category $\cO_{int}$ is semisimple with irreducible 
objects $\{V(\lm ):\lm \in P_+\},$

where $\lm \in P_{+}:=\{ \lambda \in \mathfrak{h} ^* : \lambda (P^{\vee}) \subset \bbN \}.$ 

A weight module, $M$, satisfying (i) above is said to be a \textbf{highest weight module of highest weight} $\lm$ if there exists a $u \in M$ such that
\begin{enumerate}
\item[(i)]$e_i u=0 \mbox{ for all } i \in I,$
\item[(ii)]$q^h u=q^{\lm(h)}u$ for all $h \in P^\vee, $ and
\item[(iii)]$M=\cU _q(\mathfrak g)u.$
\end{enumerate}

Every highest weight module in the category $\cO_{int}$ is isomorphic to $V(\lm)$ for some $\lm \in P_{+}.$

For $k \in I$, define $\Lambda _k \in  P_{+}$ by $\Lambda_k (h_j) = \delta _{kj},$ and $\Lambda_k(d)=0.$
  In the following sections we will define the Fock space representation  for $\Lambda_k.$  The corresponding module belongs to the category $\cO_{int}$  and contains $V(\Lambda_k).$ We will use this to describe $B(\Lambda_k)$ by coloured  Young diagrams, where $(L(\Lambda_k),B(\Lambda_k))$ denotes the (lower or upper) crystal base of $V(\Lambda_k)$ (see \cite{Ka91}).


\section*{The Fock space representations for $\mathfrak{g}$ of type $C_n^{(1)}$}
\setcounter{section}{2}
\setcounter{subsection}{0}

\label{sec2.0}

Here we modify the definitions in \cite{JMMO} to define the Fock space representations of 
$\cU _q(\mathfrak g). $

\begin{defn}
A \textbf{Young diagram }\index{Young diagram} $Y$ \textbf{of charge $k,$} for $ k \in I,$ is a 
sequence $\{y_{l}\}_{l \in \bbN}$ such that
\begin{enumerate}
    \item [(i)]
    $y_{l} \in \bbZ,$
    \item [(ii)]
    $y_{l} \leq y_{l+1}$ for all $l \in \bbN$, and 
    \item [(iii)]
    $y_{l}=i$ for all $l>>0.$
\end{enumerate}
The empty Young diagram of charge $k$ will be denoted by $\phi_k,$ i.e. $\phi_k = (k, k, \ldots ).$

\noindent
Define 
\[\cY(\Lambda_k):=\{\mathbf{Y}: \my \mbox{ is a Young diagram of charge } k \},\]
and the \textbf{Fock space} of weight $\Lambda_k$ to be
\[
\cF(\Lambda_i)=\bigoplus _{\mathbf{Y} \in \cY} \bbQ \mathbf{Y}.
\]

\end{defn}

We colour the $x$-$y$ plane as follows: For $l$ and $l^{\prime} \in \bbZ$, the ``box'' $\{(x,y) : 
 \ l<x \leq l+1, \ l^{\prime} -1< y \leq 
l^{\prime}  \}$ is coloured $i$ where $i \in I$ and 
$l+l^{\prime} \equiv \pm i \bmod 2n $ (see Figure~\ref{fig1}).  Then the diagram $Y=\{y_{l}\}_{l \geq 
0}$ is represented in the coloured $x$-$y$ plane by the coloured 
region defined by $\{(x,y) : \ l \leq x \leq l+1, \ 0 \geq y \geq y_{l}
\mbox{ for some } l \in \bbN \}.$


 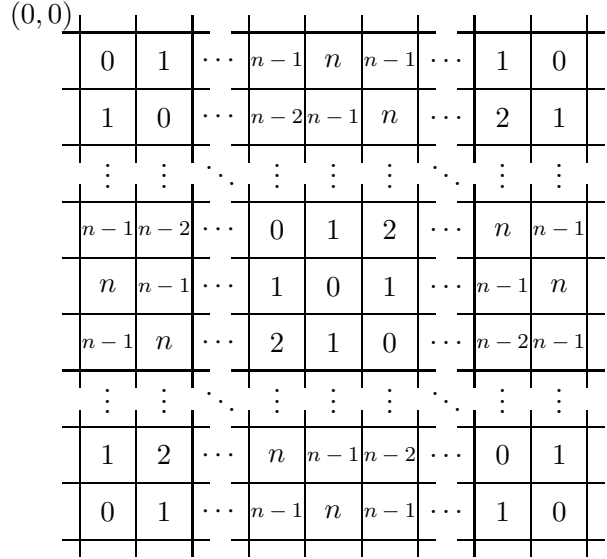
\begin{figure}
\small

\setlength{\unitlength}{0.75cm}
\begin{picture}(11,11)(-1.5,-10)
    \put(-0.3,0){\line(1,0){2.6}}
    \put(2.7,0){\line(1,0){3.6}}
    \put(6.7,0){\line(1,0){2.6}}
    
    \put(-0.3,-1){\line(1,0){2.6}}
    \put(2.7,-1){\line(1,0){3.6}}
    \put(6.7,-1){\line(1,0){2.6}}

    \put(-0.3,-2){\line(1,0){2.6}}
    \put(2.7,-2){\line(1,0){3.6}}
    \put(6.7,-2){\line(1,0){2.6}}

    \put(-0.3,-3){\line(1,0){2.6}}
    \put(2.7,-3){\line(1,0){3.6}}
    \put(6.7,-3){\line(1,0){2.6}}

    \put(-0.3,-4){\line(1,0){2.6}}
    \put(2.7,-4){\line(1,0){3.6}}
    \put(6.7,-4){\line(1,0){2.6}}

    \put(-0.3,-5){\line(1,0){2.6}}
    \put(2.7,-5){\line(1,0){3.6}}
    \put(6.7,-5){\line(1,0){2.6}}

    \put(-0.3,-6){\line(1,0){2.6}}
    \put(2.7,-6){\line(1,0){3.6}}
    \put(6.7,-6){\line(1,0){2.6}}

    \put(-0.3,-7){\line(1,0){2.6}}
    \put(2.7,-7){\line(1,0){3.6}}
    \put(6.7,-7){\line(1,0){2.6}}

    \put(-0.3,-8){\line(1,0){2.6}}
    \put(2.7,-8){\line(1,0){3.6}}
    \put(6.7,-8){\line(1,0){2.6}}

    \put(-0.3,-9){\line(1,0){2.6}}
    \put(2.7,-9){\line(1,0){3.6}}
    \put(6.7,-9){\line(1,0){2.6}}

    \put(0,0.3){\line(0,-1){2.6}}
    \put(0,-2.7){\line(0,-1){3.6}}
    \put(0,-6.7){\line(0,-1){2.6}}
    
    \put(1,0.3){\line(0,-1){2.6}}
    \put(1,-2.7){\line(0,-1){3.6}}
    \put(1,-6.7){\line(0,-1){2.6}}

    \put(2,0.3){\line(0,-1){2.6}}
    \put(2,-2.7){\line(0,-1){3.6}}
    \put(2,-6.7){\line(0,-1){2.6}}

    \put(3,0.3){\line(0,-1){2.6}}
    \put(3,-2.7){\line(0,-1){3.6}}
    \put(3,-6.7){\line(0,-1){2.6}}

    \put(4,0.3){\line(0,-1){2.6}}
    \put(4,-2.7){\line(0,-1){3.6}}
    \put(4,-6.7){\line(0,-1){2.6}}

    \put(5,0.3){\line(0,-1){2.6}}
    \put(5,-2.7){\line(0,-1){3.6}}
    \put(5,-6.7){\line(0,-1){2.6}}

    \put(6,0.3){\line(0,-1){2.6}}
    \put(6,-2.7){\line(0,-1){3.6}}
    \put(6,-6.7){\line(0,-1){2.6}}

    \put(7,0.3){\line(0,-1){2.6}}
    \put(7,-2.7){\line(0,-1){3.6}}
    \put(7,-6.7){\line(0,-1){2.6}}

    \put(8,0.3){\line(0,-1){2.6}}
    \put(8,-2.7){\line(0,-1){3.6}}
    \put(8,-6.7){\line(0,-1){2.6}}

    \put(9,0.3){\line(0,-1){2.6}}
    \put(9,-2.7){\line(0,-1){3.6}}
    \put(9,-6.7){\line(0,-1){2.6}}

   \put(0.5,-0.5){\makebox(0,0)[c]{$0$}}
   \put(1.5,-0.5){\makebox(0,0)[c]{$1$}}
   \put(2.5,-0.5){\makebox(0,0)[c]{$\cdots$}}
   \put(3.5,-0.5){\makebox(0,0)[c]{\tiny{$n-1$}}}
   \put(4.5,-0.5){\makebox(0,0)[c]{$n$}}
   \put(5.5,-0.5){\makebox(0,0)[c]{\tiny{$n-1$}}}
   \put(6.5,-0.5){\makebox(0,0)[c]{$\ldots$}}
   \put(7.5,-0.5){\makebox(0,0)[c]{$1$}}
   \put(8.5,-0.5){\makebox(0,0)[c]{$0$}}

   \put(0.5,-1.5){\makebox(0,0)[c]{$1$}}
   \put(1.5,-1.5){\makebox(0,0)[c]{$0$}}
   \put(2.5,-1.5){\makebox(0,0)[c]{$\cdots$}}
   \put(3.5,-1.5){\makebox(0,0)[c]{\tiny{$n-2$}}}
   \put(4.5,-1.5){\makebox(0,0)[c]{\tiny{$n-1$}}}
   \put(5.5,-1.5){\makebox(0,0)[c]{{$n$}}}
   \put(6.5,-1.5){\makebox(0,0)[c]{$\ldots$}}
   \put(7.5,-1.5){\makebox(0,0)[c]{$2$}}
   \put(8.5,-1.5){\makebox(0,0)[c]{$1$}}

   \put(0.5,-2.4){\makebox(0,0)[c]{$\vdots$}}
   \put(1.5,-2.4){\makebox(0,0)[c]{$\vdots$}}
   \put(2.5,-2.4){\makebox(0,0)[c]{$\ddots$}}
   \put(3.5,-2.4){\makebox(0,0)[c]{$\vdots$}}
   \put(4.5,-2.4){\makebox(0,0)[c]{$\vdots$}}
   \put(5.5,-2.4){\makebox(0,0)[c]{$\vdots$}}
   \put(6.5,-2.4){\makebox(0,0)[c]{$\ddots$}}
   \put(7.5,-2.4){\makebox(0,0)[c]{$\vdots$}}
   \put(8.5,-2.4){\makebox(0,0)[c]{$\vdots$}}

   \put(0.5,-3.5){\makebox(0,0)[c]{\tiny{$n-1$}}}
   \put(1.5,-3.5){\makebox(0,0)[c]{\tiny{$n-2$}}}
   \put(2.5,-3.5){\makebox(0,0)[c]{$\cdots$}}
   \put(3.5,-3.5){\makebox(0,0)[c]{0}}
   \put(4.5,-3.5){\makebox(0,0)[c]{$1$}}
   \put(5.5,-3.5){\makebox(0,0)[c]{$2$}}
   \put(6.5,-3.5){\makebox(0,0)[c]{$\ldots$}}
   \put(7.5,-3.5){\makebox(0,0)[c]{$n$}}
   \put(8.5,-3.5){\makebox(0,0)[c]{{\tiny$n-1$}}}

   \put(0.5,-4.5){\makebox(0,0)[c]{$n$}}
   \put(1.5,-4.5){\makebox(0,0)[c]{\tiny{$n-1$}}}
   \put(2.5,-4.5){\makebox(0,0)[c]{$\cdots$}}
   \put(3.5,-4.5){\makebox(0,0)[c]{1}}
   \put(4.5,-4.5){\makebox(0,0)[c]{$0$}}
   \put(5.5,-4.5){\makebox(0,0)[c]{$1$}}
   \put(6.5,-4.5){\makebox(0,0)[c]{$\ldots$}}
   \put(7.5,-4.5){\makebox(0,0)[c]{{\tiny$n-1$}}}
   \put(8.5,-4.5){\makebox(0,0)[c]{$n$}}

   \put(0.5,-5.5){\makebox(0,0)[c]{\tiny{$n-1$}}}
   \put(1.5,-5.5){\makebox(0,0)[c]{$n$}}
   \put(2.5,-5.5){\makebox(0,0)[c]{$\cdots$}}
   \put(3.5,-5.5){\makebox(0,0)[c]{2}}
   \put(4.5,-5.5){\makebox(0,0)[c]{$1$}}
   \put(5.5,-5.5){\makebox(0,0)[c]{$0$}}
   \put(6.5,-5.5){\makebox(0,0)[c]{$\ldots$}}
   \put(7.5,-5.5){\makebox(0,0)[c]{{\tiny$n-2$}}}
   \put(8.5,-5.5){\makebox(0,0)[c]{{\tiny$n-1$}}}

   \put(0.5,-6.4){\makebox(0,0)[c]{$\vdots$}}
   \put(1.5,-6.4){\makebox(0,0)[c]{$\vdots$}}
   \put(2.5,-6.4){\makebox(0,0)[c]{$\ddots$}}
   \put(3.5,-6.4){\makebox(0,0)[c]{$\vdots$}}
   \put(4.5,-6.4){\makebox(0,0)[c]{$\vdots$}}
   \put(5.5,-6.4){\makebox(0,0)[c]{$\vdots$}}
   \put(6.5,-6.4){\makebox(0,0)[c]{$\ddots$}}
   \put(7.5,-6.4){\makebox(0,0)[c]{$\vdots$}}
   \put(8.5,-6.4){\makebox(0,0)[c]{$\vdots$}}

   \put(0.5,-7.5){\makebox(0,0)[c]{$1$}}
   \put(1.5,-7.5){\makebox(0,0)[c]{$2$}}
   \put(2.5,-7.5){\makebox(0,0)[c]{$\cdots$}}
   \put(3.5,-7.5){\makebox(0,0)[c]{$n$}}
   \put(4.5,-7.5){\makebox(0,0)[c]{\tiny{$n-1$}}}
   \put(5.5,-7.5){\makebox(0,0)[c]{\tiny{$n-2$}}}
   \put(6.5,-7.5){\makebox(0,0)[c]{$\ldots$}}
   \put(7.5,-7.5){\makebox(0,0)[c]{$0$}}
   \put(8.5,-7.5){\makebox(0,0)[c]{$1$}}

   \put(0.5,-8.5){\makebox(0,0)[c]{$0$}}
   \put(1.5,-8.5){\makebox(0,0)[c]{$1$}}
   \put(2.5,-8.5){\makebox(0,0)[c]{$\cdots$}}
   \put(3.5,-8.5){\makebox(0,0)[c]{\tiny{$n-1$}}}
   \put(4.5,-8.5){\makebox(0,0)[c]{$n$}}
   \put(5.5,-8.5){\makebox(0,0)[c]{\tiny{$n-1$}}}
   \put(6.5,-8.5){\makebox(0,0)[c]{$\ldots$}}
   \put(7.5,-8.5){\makebox(0,0)[c]{$1$}}
   \put(8.5,-8.5){\makebox(0,0)[c]{$0$}}

    \put(-0.1,0.3){\makebox(0,0)[r]{$(0,0)$}}

    \normalsize
\end{picture}
\caption{The colouring of the x-y plane.}
\label{fig1}
\end{figure}


Let $\mathbf{Y}=\{y_l\}_{l \in \bbN}  \in \cY(\Lambda_k).$   If $y_{l} \not = 
y_{l+1}$ for some $l \in \bbN$, 
$\mathbf{Y}$ is said to have a \textbf{concave} (\textbf{convex})
\index{concave corner} \index{convex corner}
\textbf{corner} at site $(l+1, y_{l+1})$ ($(l+1, y_{l})$, 
resp.).  Also , $\mathbf{Y}$ is said to 
have a \textbf{concave corner} at site $(0, y_{l})$.  
For $i \in I,$ a corner at site $(l, y)$ is called an $i-$\textbf{coloured corner}
\index{ $i-$coloured corner }
if  $l+y \equiv \pm i \bmod 2n$.
\smallskip


\begin{sex}\label{sec3.2.1}
    Let $n=2 \mas$
    $\mathbf{Y}=(-4,-2,-2,-1,-1,0,0,\ldots).$
\[
\begin{array}{ll}
\begin{array}{l}
 \my= \\ \\ \\ \\ \\ \\ 
\end{array}

\begin{Young}
0&1&2&1&0 \cr
1&0&1 \cr
2 \cr
1 \cr
\end{Young}
\end{array}
\]

   {\my} has a $0$-coloured concave corner 
    at sites $(0,-4),$   a $0$-coloured convex corner at site $(5,-1),$  $1$-coloured concave corners at sites 
    $(0,5) \mas (1,-2),$  $1$-coloured convex corners at sites $(3,-2), \mas (1,-4)$, and a 
    $2$-coloured concave corner at sites $(3, -1).$

\end{sex}	


\subsection{} \label{sec3.4}
				
We now define an action of $\cU_q(\mathfrak {g})$ on $\cF(\Lambda_k).$

For $(l,y) \in  \bbN \times \bbZ$ define  linear maps $E_{(l,y)},\ F_{(l,y)},\ T_{(l,y)}^{\pm}:\cF(\Lambda_k) \to \cF(\Lambda_k)$ as follows, for $\mathbf{Y} \in \cY(\Lambda_k),$

if $\my$ has a convex corner at site $(l,y), E_{(l,y)} (\mathbf{Y})$ is the same as $\mathbf{Y}$ with this corner removed; otherwise, $E_{(l,y)} (\mathbf{Y})=0,$

if $\bf{Y}$ has a concave corner at site $(l,y), F_{(l,y)} (\mathbf{Y})$ is the same as $\mathbf{Y}$ with a corner added at site (l+1,y-1) ; otherwise, $F_{(l,y)} (\mathbf{Y})=0,$

\[
T_{(l,y)}^{\pm}(\mathbf{Y})=
\begin{cases}
    q_i^{\pm} \mathbf{Y} & \mbox{ if } \mathbf{Y} \mbox{ has a concave corner at site } $(l,y)$,
     \\
    q_i^{\mp} \mathbf{Y} & \mbox{ if } \mathbf{Y} \mbox{ has a convex corner at site } $(l,y)$,
     \\
    \my & \mbox{otherwise,}
\end{cases}\]
\noindent
where $i \in I$ such that $l+y \equiv \pm i\bmod 2n.$

Define the order $>$ on $\bbN \times \bbZ$ as follows:
$ (l, y) > (l^{\prime}, y^{\prime})  \mbox{ iff } l+y > l^{\prime} + y^{\prime}. $
 (See Figure \ref{fig2} where a point lying on a diagonal line labeled by $a$ is greater than a point lying on a diagonal line labeled by $b$ if $a>b.$)

 \begin{figure}
\small

\setlength{\unitlength}{.75cm}
\begin{picture}(6,6)(4,-9)

    \put(3,-3){\line(1,-1){5.6}}
    \put(3,-4){\line(1,-1){4.6}}
    \put(3,-5){\line(1,-1){3.6}}
    \put(3,-6){\line(1,-1){2.6}}
    \put(3,-7){\line(1,-1){1.6}}
    \put(3,-8){\line(1,-1){0.6}}

    \put(3.7,-2.7){\line(1,-1){4.9}}
    \put(4.7,-2.7){\line(1,-1){3.9}}
    \put(5.7,-2.7){\line(1,-1){2.9}}
    \put(6.7,-2.7){\line(1,-1){1.9}}
    \put(7.7,-2.7){\line(1,-1){0.9}}

    \put(2.7,-3){\line(1,0){5.6}}
    \put(2.7,-4){\line(1,0){5.6}}
    \put(2.7,-5){\line(1,0){5.6}}
    \put(2.7,-6){\line(1,0){5.6}}
    \put(2.7,-7){\line(1,0){5.6}}
    \put(2.7,-8){\line(1,0){5.6}}

    \put(3,-2.7){\line(0,-1){5.6}}
    \put(4,-2.7){\line(0,-1){5.6}}
    \put(5,-2.7){\line(0,-1){5.6}}
    \put(6,-2.7){\line(0,-1){5.6}}
    \put(7,-2.7){\line(0,-1){5.6}}
    \put(8,-2.7){\line(0,-1){5.6}}

   \put(9,-3.7){\makebox(0,0)[c]{5}}
   \put(9,-4.7){\makebox(0,0)[c]{4}}
   \put(9,-5.7){\makebox(0,0)[c]{3}}
   \put(9,-6.7){\makebox(0,0)[c]{2}}
   \put(9,-7.7){\makebox(0,0)[c]{1}}
   \put(9,-8.7){\makebox(0,0)[c]{0}}
   \put(8,-8.7){\makebox(0,0)[c]{-1}}
   \put(7,-8.7){\makebox(0,0)[c]{-2}}
   \put(6,-8.7){\makebox(0,0)[c]{-3}}
   \put(5,-8.7){\makebox(0,0)[c]{-4}}
   \put(4,-8.7){\makebox(0,0)[c]{-5}}





    \put(2.7,-2.7){\makebox(0,0)[r]{\tiny{$(0,0)$}}}

    \normalsize
\end{picture}
\caption{The order of $\bbN \times \bbZ$.}
\label{fig2}
\end{figure}


Now define linear operators $E_i,\ F_i,\  T_i, i \in I, \mas T_d : \cF(\Lambda_k) \to \cF(\Lambda_k)$ as follows,

\[
E_i=\sum_{\stackrel{(l,y) \in  \times \bbN \times \bbZ}{l+y \equiv \pm i \bmod 2n}} \big (\prod_{\stackrel{(l^{\prime},y^{\prime})>(l,y)}{l^{\prime}+y^{\prime} \equiv \pm i \bmod 2n}}T_{(l^{\prime},y^{\prime})}^+ \big ) E_{(l,y)},
\]

\[
F_i=\sum_{\stackrel{(l,y) \in  \times \bbN \times \bbZ}{l+y \equiv \pm i \bmod 2n}} \big (\prod_{\stackrel{(l^{\prime},y^{\prime})<(l,y)}{l^{\prime}+y^{\prime} \equiv \pm i \bmod 2n}}T_{(l^{\prime},y^{\prime})}^- \big ) F_{(l,y)},
\]

\[
T_i^{\pm}=\prod_{\stackrel{(l,y) \in  \times \bbN \times \bbZ}{l+y \equiv \pm i \bmod 2n}} T_{(l,y)}^{\pm}, \mas
\]

\[
\mbox{for } \mathbf{Y} \in \cY(\lambda),
T_d (\mathbf{Y})=q^{\mbox{-(the number of 0-coloured boxes in Y)}}
\]

\bigskip

\begin{thm}
 The vector space $\cF(\Lambda_k),$  where $k \in I,$ is a $\, \cU _q (\mathfrak{g})$-module where the action of the generators $e_i,\ f_i,\ q^{h_i}$ and $q^d$ is  given by that of $E_i,\ F_i,\ T_i$ and $T_d,$ respectively.
\end{thm}

\begin{pf}
We have to show that $E_i,\ F_i,\ T_i$ and $T_d$ satisfy the defining relations of $e_i,\ f_i,\ q^{h_i}$ and $q^d.$

(1) is clear.
\bigskip

(2)For $i, j \in I,\ T_jE_i=q^{a_{ji}} E_i T_j$ will follow from $T_jE_{(l,y)}=q^{a_{ji}} E_{(l,y)} T_j$ for 
$(l,y) \in \bbN \times \bbZ$ with $l+y \equiv \pm i \bmod 2n.$

Note that, for $\mathbf{Y} \in \cY(\Lambda_k) ,T_j(\mathbf{Y})=q^{ccj(\mathbf{Y}) - cxj(\mathbf{Y})}, $ where
$ccj(\mathbf{Y}):=\#$ of concave j-coloured corners in $ \mathbf{Y}$ and 
$cxj(\mathbf{Y}):=\#$of convex j-coloured corners in  $\mathbf{Y}.$

Assume that $E_{(l,y)}(\mathbf{Y}) \not = 0.$

If $i \mas j$ are not adjacent nodes in the Dynkin diagram and $i \not = j$, then $ccj(\mathbf{Y}) - cxj(\mathbf{Y})=ccj(E_{(l,y)}(\mathbf{Y})) - cxj(E_{(l,y)}(\mathbf{Y})).$

If $i \mas j$ are adjacent nodes in the Dynkin diagram and $i \not = 0 
\mbox{ or } n$, then $E_i(\mathbf{Y})$ has  one less concave $j$-coloured corner or one more convex $j$-coloured corner than $\mathbf{Y}$ (Figure \ref{fig3}  shows all the possible cases.)  In either case $ccj(E_{(l,y)}(\mathbf{Y})) - cxj(E_{(l,y)}(\mathbf{Y}))=ccj(\mathbf{Y}) - cxj(\mathbf{Y})-1.$  

 \begin{figure}
\small

\setlength{\unitlength}{.75cm}
\begin{picture}(10,3)(2,-1)

    \put(1,1){\line(1,0){2}}
    \put(1,0){\line(1,0){1}}
    \put(1,1){\line(0,-1){1}}
    \put(2,1){\line(0,-1){1}}
    \put(1.5,.5){\makebox(0,0)[c]{i}}
    \put(2.5,.5){\makebox(0,0)[c]{j}}
   
    \put(5,1){\line(1,0){1}}
    \put(5,0){\line(1,0){1}}
    \put(5,1){\line(0,-1){1}}
    \put(6,0){\line(0,1){2}}
    \put(5.5,0.5){\makebox(0,0)[c]{i}}
    \put(5.5,1.5){\makebox(0,0)[c]{j}}

    \put(8,1){\line(1,0){1}}
    \put(8,0){\line(1,0){1}}
    \put(8,1){\line(0,-1){2}}
    \put(9,0){\line(0,1){1}}
    \put(8.5,0.5){\makebox(0,0)[c]{i}}
    \put(8.5,-.5){\makebox(0,0)[c]{j}}

    \put(12,1){\line(1,0){1}}
    \put(11,0){\line(1,0){2}}
    \put(12,1){\line(0,-1){1}}
    \put(13,1){\line(0,-1){1}}
    \put(12.5,0.5){\makebox(0,0)[c]{i}}
    \put(11.5,.5){\makebox(0,0)[c]{j}}
    \normalsize
\end{picture}
\caption{}
\label{fig3}
\end{figure}


If $i=0 \mas j=1$ (or $i=n$ and $j=n-1$), then $E_i(\mathbf{Y})$ has  two less concave j-coloured corners or  one less concave $j$-coloured corner and one more convex $j$-coloured corner or  two more convex $j$-coloured corners than $\mathbf{Y}$ (Figure \ref{fig4}  shows all the possible cases.)  In all case $ccj(E_{(l,y)}(\mathbf{Y})) - cxj(E_{(l,y)}(\mathbf{Y}))=ccj(\mathbf{Y}) - cxj(\mathbf{Y})-2.$

 \begin{figure}
\small

\setlength{\unitlength}{.75cm}
\begin{picture}(10,3)(2,-1)

    \put(0,1){\line(1,0){2}}
    \put(0,0){\line(1,0){1}}
    \put(0,1){\line(0,-1){2}}
    \put(1,1){\line(0,-1){1}}
    \put(0.5,.5){\makebox(0,0)[c]{i}}
    \put(1.5,.5){\makebox(0,0)[c]{j}}
    \put(0.5,-.5){\makebox(0,0)[c]{j}}
   
    \put(5,1){\line(1,0){1}}
    \put(4,0){\line(1,0){2}}
    \put(5,1){\line(0,-1){1}}
    \put(6,0){\line(0,1){2}}
    \put(5.5,0.5){\makebox(0,0)[c]{i}}
    \put(5.5,1.5){\makebox(0,0)[c]{j}}
    \put(4.5,0.5){\makebox(0,0)[c]{j}}

    \put(8,1){\line(1,0){1}}
    \put(8,0){\line(1,0){1}}
    \put(8,1){\line(0,-1){2}}
    \put(9,0){\line(0,1){2}}
    \put(8.5,0.5){\makebox(0,0)[c]{i}}
    \put(8.5,-.5){\makebox(0,0)[c]{j}}
    \put(8.5,1.5){\makebox(0,0)[c]{j}}

    \put(12,1){\line(1,0){1}}
    \put(11,0){\line(1,0){2}}
    \put(12,1){\line(0,-1){1}}
    \put(13,2){\line(0,-1){2}}
    \put(12.5,0.5){\makebox(0,0)[c]{i}}
    \put(11.5,.5){\makebox(0,0)[c]{j}}
    \put(12.5,1.5){\makebox(0,0)[c]{j}}
    \normalsize
\end{picture}
\caption{}
\label{fig4}
\end{figure}


If $i=j$, $E_i(\mathbf{Y})$ has one less convex $i$-coloured corner and one more concave i-coloured corner than $\mathbf{Y}$.  Hence $ccj(E_{(l,y)}(\mathbf{Y})) - cxj(E_{(l,y)}(\mathbf{Y}))=ccj(\mathbf{Y}) - cxj(\mathbf{Y})+2.$

We now show that $T_d E_i =q^{\delta_{0,i}} E_i T_d.$  If $i \not =0,\ T_d \mas E_i$ commute.  If $i=0$ and $\mathbf{Y} \in \cY(\Lambda_k)$ then either $E_0(\mathbf{Y})=0$ or 

$T_d E_0 (\mathbf{Y})=q^{-\mbox{ the } \# \mbox{ of } 0-\mbox{coloured boxes in }\mathbf{Y}+1}E_0(\mathbf{Y}). $ In either case the result follows.

The proof of (3) is similar.
\bigskip

(4)  For $\my \in \cY(\Lambda_k),\ i \in I \mas (l,y) \in \bbN \times \bbZ,$ define 
\begin{align*}
a(i,l,y,\my)=\#\{(l^\prime,y^\prime) \ : &(l^\prime,y^\prime) \mbox{ is a  concave }i-\mbox{coloured corner in }\my \mas \\
&(l^\prime,y^\prime)> (l,y)\} \\
-\#\{(l^\prime,y^\prime) \ : &(l^\prime,y^\prime) \mbox{ is a  convex }i-\mbox{coloured corner in }\my \mas \\
&(l^\prime,y^\prime)> (l,y)\} \mas
\end{align*}
\begin{align*}
b(i,l,y,\my)=\#\{(l^\prime,y^\prime) \ : &(l^\prime,y^\prime) \mbox{ is a  convex }i-\mbox{coloured corner in }\my \mas \\
&(l^\prime,y^\prime)< (l,y)\} \\
-\#\{(l^\prime,y^\prime) \ : &(l^\prime,y^\prime) \mbox{ is a  concave }i-\mbox{coloured corner in }\my \mas \\
&(l^\prime,y^\prime)< (l,y)\}.
\end{align*}
Then \[
E_i(\my) =\sum_{\stackrel{(l,y) \in \bbN \times \bbZ}{l+y \equiv \pm i \bmod 2n}} 
q_i^{a(i,l,y,\my)}E_{(l,y)}(\my),
\]
\[
F_i(\my) =\sum_{\stackrel{(l,y) \in \bbN \times \bbZ}{l+y \equiv \pm i \bmod 2n}} 
q_i^{b(i,l,y,\my)}F_{(l,y)}(\my),
\]
\[
E_i(F_j(\my)) =\sum_{\stackrel{(l,y) \in}{\stackrel{ \bbN \times \bbZ}{l+y \equiv  \pm i \bmod 2n}}}
\sum_{\stackrel{(l_1,y_1) \in}{\stackrel { \bbN \times \bbZ}{l_1+y_1 \equiv  \pm j \bmod 2n}}}
q_i^{a(i,l,y,F_{(l_1,y_1)}(\my))}q_j^{b(j,l_1,y_1,\my)}E_{(l,y)}(F_{(l_1,y_1)}(\my)),
\]
\[
F_j(E_i(\my)) =\sum_{\stackrel{(l,y) \in}{\stackrel{ \bbN \times \bbZ}{l+y \equiv  \pm i \bmod 2n}}}
\sum_{\stackrel{(l_1,y_1) \in}{\stackrel { \bbN \times \bbZ}{l_1+y_1 \equiv  \pm j \bmod 2n}}}
q_i^{a(i,l,y,(\my))}q_j^{b(j,l_1,y_1,E_{(l_1,y_1)}(\my))}F_{(l_1,y_1)}(E_{(l,y)}(\my)),
\]

\noindent
In what follows, $(l,y), (l_1,y_1) \in \bbN \times \bbZ$ with $l+y \equiv  \pm i \bmod 2n ,
\ l_1+y_1 \equiv  \pm j \bmod 2n.$
Note that $E_{(l,y)}F_{(l_1,y_1)}= F_{(l_1,y_1)} E_{(l,y)}$ unless $(l_1,y_1)=(l-1,y+1).$
We will assume that $E_{(l,y)}F_{(l_1,y_1)}(\my) \not = 0$ or $F_{(l_1,y_1)} E_{(l,y)}(\my) \not = 0.$

If $i \mas j$ are not adjacent in the Dynkin diagram and $i \not = j,$ or if $(l_1,y_1) < (l,y)$
\begin{align*}
&a(i,l,y,F_{(l_1,y_1)}(\my))=a(i,l,y,(\my)),\\
&b(j,l_1,y_1,\my)=b(j,l_1,y_1,E_{(l,y)}(\my)).
\end{align*}

If $i \mas j$ are adjacent and $(l_1,y_1) > (l,y),$ we have three cases.

\textbf{Case 1.} $j=0 \mas i=1$ (or $j=n \mas i=n-1).$ 
\begin{align*}
&a(i,l,y,F_{(l_1,y_1)}(\my))=a(i,l,y,(\my))+2,\\
&b(j,l_1,y_1,\my)=b(j,l_1,y_1,E_{(l,y)}(\my))-1.
\end{align*}

\textbf{Case 2.} $i=0 \mas j=1$ (or $i=n \mas j=n-1).$ 
\begin{align*}
&a(i,l,y,F_{(l_1,y_1)}(\my))=a(i,l,y,(\my))+1,\\
&b(j,l_1,y_1,\my)=b(j,l_1,y_1,E_{(l,y)}(\my))-2.
\end{align*}

\textbf{Case 3.} $j \not =0,n,\  i \not = 0,n \mas i \not = j.$  
\begin{align*}
&a(i,l,y,F_{(l_1,y_1)}(\my))=a(i,l,y,(\my))+1,\\
&b(j,l_1,y_1,\my)=b(j,l_1,y_1,E_{(l,y)}(\my))-1.
\end{align*}

Hence if $i \not = j$, $E_iF_j=F_jE_i.$

\bigskip
If $i =j \mas (l_1,y_1) > (l,y),$
\begin{align*}
&a(i,l,y,F_{(l_1,y_1)}(\my))=a(i,l,y,(\my))-2,\\
&b(j,l_1,y_1,\my)=b(j,l_1,y_1,E_{(l,y)}(\my))+2.
\end{align*}

So the only terms in $(E_iF_i-F_iE_i)(\my)$ which give  non-zero contributions to the sum are those terms where
$ (l_1,y_1)=(l-1,y+1),$ and $E_{(l,y)}(F_{(l_1,y_1)}(\my)) \not = 0$ or $F_{(l_1,y_1)}(E_{(l,y)}(\my))\not =0.$

A concave corner in $\my$ at site $(l,y)$ will contribute $q_i ^{a(i,l,y,F_{(l,y)}(\my))b(i,l,y,\my)}\my$ to $(E_iF_i-F_iE_i)(\my)$ and
a convex corner in $\my$ at site $(l,y)$ a $\ -q_i ^{a(i,l,y,\my)b(i,l,y,E_{(l,y)}(\my))}\my.$

Let $(l_1,y_1) >(l_2,y_2)>\cdots >(r_l,k_l,y_l)$ be the sites of the $i$-coloured corners of $\my,$ $\sigma$ be the $i$-signature of $\my$ (see \ref{sec3.2})
and $J(\sigma)$ be defined as in section \ref{sec3.2}.  Then
the contribution to $(E_iF_i-F_iE_i)(\my)$ of the corners in $\my$ corresponding to $\{1,2,...,l\} \backslash J(\sigma)$ cancel out.  Let $a=\# \mbox { of } 1's \mbox { in } J(\sigma)$ and $b=\# \mbox { of } 0's \mbox { in } J(\sigma).$   then 
\[
(E_iF_i-F_iE_i)(\my)=\frac{q_i^{b-a}-q_i^{-(b-a)}}{q_i-q_i^{-1}}\my=\frac{T_i^+-T_i^-}{q_i-q_i^{-1}}\my
\]

\bigskip

(5) If $a_{ij}=0$ i.e. $i \mas j$ are not adjacent, $E_iE_j=E_jE_i.$

If $a_{ij}=-1,$ and $\my \in \cY(\Lambda_k),$

\begin{align*}
E_i^2E_j (\my)&=\sum_\mx a(\mx) \mx,\\
E_iE_jE_i (\my)&=\sum_\mx b(\mx) \mx,\\
E_j E_i^2(\my)&=\sum_\mx c(\mx) \mx,
\end{align*}

\noindent
where $a(\mx), b(\mx) \mas c(\mx) \in \bbQ(q)$ and the sums run over all $\mx \in \cY(\Lambda_k)$ which are obtained from $\my \in \cY(\Lambda_k)$ by removing two $i$-coloured boxes and one $j$-coloured box. 

Let $(l_1,y_1), (l_2,y_2) \mas (l_3,y_3)$ be the co-ordinates of the bottom right corners of two distinct  $i$-coloured boxes and one  $j$-coloured box of $\my$ (if they exist), resp.,  and let $\mx$ be obtained from $\my$ by removing these three boxes.  Assume $\mx \in \cY(\Lambda_k).$  

\textbf{Case1.} If $(l_1,y_1),(l_2,y_2) \mas (l_3,y_3)$ are sites of two $i$-coloured convex corners and one $j$- coloured convex corner of $\my,$ we consider three cases.

\textbf{Case 1 (a).} $(l_1,y_1),(l_2,y_2)>(l_3,y_3).$
\begin{align*}
a(\mx)&=q^m([2]_{q_i}),\\
b(\mx)&=q^mq_j^{-a_{ji}}([2]_{q_i})=q^mq_i^{-1}([2]_{q_i}), \mas\\
c(\mx)&=q^mq_j^{-2a_{ji}}([2]_{q_i})=q^mq_i^{-2}([2]_{q_i}),\mbox{ for some } m \in \bbZ .\\
\end{align*}

\textbf{Case 1 (b).} $(l_1,y_1)>(l_3,y_3)>(l_2,y_2).$
\begin{align*}
a(\mx)&=q^mq_i^{-1}([2]_{q_i}),\\
b(\mx)&=q^m(q_j^{-a_{ji}}q_i^{-1}q_i+q_i^{-1})=q^m(2q_i^{-1}),\\
c(\mx)&=q^mq_j^{-a_{ji}}([2]_{q_i})=q^mq_i^{-1}([2]_{q_i}),\mbox{ for some } m \in \bbZ .\\
\end{align*}

\textbf{Case 1 (c).} $(l_3,y_3)>(l_1,y_1),(l_2,y_2).$
\begin{align*}
a(\mx)&=q^mq_i^{-2}([2]_{q_i}),\\
b(\mx)&=q^mq_i^{-1}([2]_{q_i}),\\
c(\mx)&=q^m([2]_{q_i}),\mbox{ for some } m \in \bbZ .\\
\end{align*}

\textbf{Case2.} If the $j$-coloured box corresponding to $(l_3,y_3)$ is \textbf{hidden} by (i.e. if it is immediately to the left or above of)  the $i$-coloured box corresponding to $(l_2,y_2),$  we consider two cases.

\textbf{Case 2 (a).} $(l_1,y_1) > (l_2,y_2), (l_3,y_3).$
\begin{align*}
a(\mx)&=0,\\
b(\mx)&=q^mq_i^{-1}, \mas\\
c(\mx)&=q^mq_j^{-a_{ji}}([2]_{q_i})=q^mq_i^{-1}([2]_{q_i}),\mbox{ for some } m \in \bbZ .\\
\end{align*}

\textbf{Case 2 (b).} $(l_1,y_1)<(l_2,y_2), (l_3,y_3).$
\begin{align*}
a(\mx)&=0,\\
b(\mx)&=q^m,\\
c(\mx)&=q^m[2]_{q_i},\mbox{ for some } m \in \bbZ .\\
\end{align*}

\textbf{Case3.} If the $i$-coloured box corresponding to $(l_2,y_2)$ is hidden by the  $j$-coloured box corresponding to $(l_3,y_3),$  we consider two cases.

\textbf{Case 3 (a).} $(l_1,y_1) > (l_2,y_2), (l_3,y_3).$
\begin{align*}
a(\mx)&=q^m[2]_{q_i},\\
b(\mx)&=q^mq_j^{-a_{ji}}q_i=1, \mas\\
c(\mx)&=0,\mbox{ for some } m \in \bbZ .\\
\end{align*}

\textbf{Case 3 (b).} $(l_1,y_1)<(l_2,y_2), (l_3,y_3).$
\begin{align*}
a(\mx)&=q^m[2]_{q_i},\\
b(\mx)&=q^m,\\
c(\mx)&=0,\mbox{ for some } m \in \bbZ .\\
\end{align*}

In all cases $a(\mx)-[2]_{q_i}b(\mx)+c(\mx)=0.$



%


%


%


\bigskip

If $a_{ij}=-2,$ and $\my \in \cY(\Lambda_k),$

\begin{align*}
E_i^3E_j (\my)&=\sum_\mx a(\mx) \mx,\\
E_i^2E_jE_i (\my)&=\sum_\mx b(\mx) \mx,\\
E_iE_jE_i ^2(\my)&=\sum_\mx c(\mx) \mx,\\
E_j E_i^3(\my)&=\sum_\mx d(\mx) \mx,
\end{align*}

\noindent
where $a(\mx), b(\mx) c(\mx) \mas d(\mx) \in \bbQ(q)$ and the sums run over all $\mx \in \cY(\Lambda_k)$ which are obtained from $\my$ by removing three $i$-coloured boxes and one  $j$-coloured box. 

Let $(l_1,y_1), (l_2,y_2), (l_3,y_3) \mas (l_4,y_4)$ be the co-ordinates of the bottom right corners of  three distinct  $i$-coloured boxes and one  $j$-coloured box of $\my$ (if they exist), resp.,  and let $\mx$ be obtained from $\my$ by removing these four boxes.  Assume $\mx \in \cY({\Lambda_k}).$ 

\textbf{Case 1.} If  $\my$ has three convex $i$-coloured corners and one convex $j$-coloured corner at sites $(l_1,y_1),(l_2,y_2),(l_3,y_3) \mas (l_4,y_4),$ 
we consider four cases.

\textbf{Case 1 (a).} $(l_1,y_1),(l_2,y_2),(l_3,y_3)>(l_4,y_4).$
\begin{align*}
a(\mx)&=q^m([3]_q[2]_q),\\
b(\mx)&=q^mq^{-2}([3]_q[2]_q),\\
c(\mx)&=q^mq^{-4}([3]_q[2]_q),\\
d(\mx)&=q^mq^{-6}([3]_q[2]_q),\mbox{ for some } m \in \bbZ .\\
\end{align*}

\textbf{Case 1 (b).} $(l_1,y_1),(l_2,y_2)>(l_4,y_4)>(l_3,y_3).$
\begin{align*}
a(\mx)&=q^mq^{-2}([3]_q[2]_q),\\
b(\mx)&=q^m(2q^{-1}+3q^{-3}+q^{-5}),\\
c(\mx)&=q^m(q^{-1}+3q^{-3}+2q^{-5}),\\
d(\mx)&=q^mq^{-4}([3]_q[2]_q),\mbox{ for some } m \in \bbZ .\\
\end{align*}

\textbf{Case 1 (c).} $(l_1,y_1)>(l_4,y_4)>(l_2,y_2),(l_3,y_3).$
\begin{align*}
a(\mx)&=q^mq^{-4}([3]_q[2]_q),\\
b(\mx)&=q^m(q^{-1}+3q^{-3}+2q^{-5}),\\
c(\mx)&=q^m(2q^{-1}+3q^{-3}+q^{-5}),\\
d(\mx)&=q^mq^{-2}([3]_q[2]_q),
\mbox{ for some } m \in \bbZ .\\
\end{align*}

\textbf{Case 1 (d).} $(l_4,y_4)>(l_1,y_1),(l_2,y_2),(l_3,y_3).$
\begin{align*}
a(\mx)&=q^mq^{-6}([3]_q[2]_q),\\
b(\mx)&=q^mq^{-4}([3]_q[2]_q),\\
c(\mx)&=q^mq^{-2}([3]_q[2]_q),\\
d(\mx)&=q^m([3]_q[2]_q),\mbox{ for some } m \in \bbZ .\\
\end{align*}

\textbf{Case 2.} If  $\my$ has three convex $i$-coloured corners at sites $(l_1,y_1),(l_2,y_2), \mas(l_3,y_3)$ and the $j$-coloured box corresponding to $(l_4,y_4)$  is hidden by that corresponding to $(l_3,y_3)$ but not by the other $i$-coloured boxes, we consider three cases.

\textbf{Case 2 (a).} $(l_1,y_1),(l_2,y_2)>(l_3,y_3),(l_4,y_4).$
\begin{align*}
a(\mx)&=0,\\
b(\mx)&=q^mq^{-2}[2]_q,\\
c(\mx)&=q^m(q^{-4}+q^{-2})[2]_q,\\
d(\mx)&=q^mq^{-4}([3]_q[2]_q),\mbox{ for some } m \in \bbZ .\\
\end{align*}

\textbf{Case 2 (b).} $(l_1,y_1)>(l_3,y_3),(l_4,y_4)>(l_2,y_2).$
\begin{align*}
a(\mx)&=0,\\
b(\mx)&=q^mq^{-1}[2]_q,\\
c(\mx)&=q^m(q^{-1}+q^{-2})[2]_q,\\
d(\mx)&=q^mq^{-2}([3]_q[2]_q),\mbox{ for some } m \in \bbZ .\\
\end{align*}

\textbf{Case 2 (c).} $(l_3,y_3),(l_4,y_4)>(l_1,y_1),(l_2,y_2),.$
\begin{align*}
a(\mx)&=0,\\
b(\mx)&=q^mq^{-2}[2]_q,\\
c(\mx)&=q^m(1+q^{-2})[2]_q,\\
d(\mx)&=q^m([3]_q[2]_q),
\mbox{ for some } m \in \bbZ .\\
\end{align*}

\textbf{Case 3.} If  $\my$ has three convex $i$-coloured corners at sites $(l_1,y_1),(l_2,y_2), \mas(l_3,y_3)$ and the $j$-coloured box corresponding to $(l_4,y_4)$  is hidden by those corresponding to $(l_2,y_2) \mas (l_3,y_3),$ then

\begin{align*}
a(\mx)&=0,\\
b(\mx)&=0,\\
c(\mx)&=q^m[2]_q\\
d(\mx)&=q^m[3]_q[2]_q,\mbox{ for some } m \in \bbZ .\\
\end{align*}

\textbf{Case 4.} If  $\my$ has two convex $i$-coloured corners at sites $(l_1,y_1) \mas (l_2,y_2)$ and one $j$-coloured corner at site $(l_4,y_4),$ and the $i$-coloured box corresponding to $(l_3,y_3)$  is hidden by that corresponding to $(l_4,y_4),$ we consider three cases. 

\textbf{Case 4 (a).} $(l_1,y_1),(l_2,y_2)>(l_3,y_3),(l_4,y_4).$
\begin{align*}
a(\mx)&=q^m[3]_q[2]_q,\\
b(\mx)&=q^m(1+q^{-2})[2]_q,\\
c(\mx)&=q^mq^{-2}[2]_q,\\
d(\mx)&=0,\mbox{ for some } m \in \bbZ .\\
\end{align*}

\textbf{Case 4 (b).} $(l_1,y_1)>(l_3,y_3),(l_4,y_4)>(l_2,y_2).$
\begin{align*}
a(\mx)&=q^mq^{-1}[3]_q[2]_q,\\
b(\mx)&=q^m(2q^{-1})[2]_q,\\
c(\mx)&=q^mq^{-1}[2]_q,\\
d(\mx)&=0,\mbox{ for some } m \in \bbZ .\\
\end{align*}

\textbf{Case 4 (c).} $(l_3,y_3),(l_4,y_4)>(l_1,y_1),(l_2,y_2),.$
\begin{align*}
a(\mx)&=q^mq^{-2}[3]_q[2]_q,\\
b(\mx)&=q^m(1+q^{-2})[2]_q,\\
c(\mx)&=q^m[2]_q,\\
d(\mx)&=0,\mbox{ for some } m \in \bbZ .\\
\end{align*}

\textbf{Case 5.} If  $\my$ has one convex $i$-coloured corners at site $(l_1,y_1),$ one convex $j$-coloured corner at site $(l_4,y_4),$  and the $i$-coloured boxes corresponding to $(l_2,y_2) \mas (l_3,y_3)$  are hidden by that corresponding to $(l_4,y_4),$ then

\begin{align*}
a(\mx)&=q^m[3]_q[2]_q,\\
b(\mx)&=q^m[2]_q,\\
c(\mx)&=0\\
d(\mx)&=0,\mbox{ for some } m \in \bbZ .\\
\end{align*}

\textbf{Case 6.} If  $\my$ has two convex $i$-coloured corners at sites $(l_1,y_1) \mas (l_2,y_2),$ and the $j$-coloured box corresponding to $(l_4,y_4)$  is hidden by that corresponding to $(l_2,y_2)$ and the $i$-coloured box corresponding to $(l_3,y_3)$ is hidden by that corresponding to $(l_4,y_4)$ then

\begin{align*}
a(\mx)&=0,\\
b(\mx)&=q^m[2]_q,\\
c(\mx)&=q^m[2]_q,\\
d(\mx)&=0,\mbox{ for some } m \in \bbZ .\\
\end{align*}

In all cases $a(\mx)-[3]_{q_i}b(\mx)+[3]_{q_i}c(\mx)-d(\mx)=0.$

\bigskip

(6) is similar to (5).
\end{pf}

\begin{lem}
The $\cU_q(\mathfrak{g})$-module $\cF(\Lambda_k)$ belongs to the category $\cO_{int}.$
\end{lem}

\begin{pf}
Note that $q^h \phi_k = q^{\Lambda_k(h)}\phi_k$ for all $h \in P^\vee.$   By induction on the number of boxes of an element $\my \in \cY(\Lambda_k),$  we can show that $q^h \my=q^{\mu(h)}\my$ for all $h \in P^\vee,$  where $ \mu=\Lambda_k-\sum_{i=0}^n k_i \alpha_i$ and $k_i=\#$ of $i$-coloured boxes in $\my.$

\end{pf}

\begin{cor}
$M(\Lambda_k):=\cU_q(\mathfrak{g})\phi_k$ is the irreducible integrable highest weight module of highest weight $\Lambda_k.$
\end{cor}

\section*{The Crystal Base for  $\cF(\Lambda_k)$  and the crystal $B(\Lambda_k)$}

\setcounter{section}{3}
\setcounter{subsection}{0}


The proofs of the following two theorems are as in \cite{JMMO} and \cite{MM90}.

\begin{thm}
Let $A=\{\frac{f(q)}{g(q)}:\ f(q),g(q) \in \bbQ[q] \mas g(0) \not =0\}, k \in I, \ L(\cF(\Lambda_k))=\sum_{\my \in \cF(\Lambda_k)}A \my$ and $B(\cF(\Lambda_k))=\cY(\Lambda_k),$ where we identify $\my +qL(\cF(\Lambda_k))$ with $ \my.$
Then $(L(\cF(\Lambda_k)),B(\cF(\Lambda_k)))$ is an (upper) crystal base for the integrable $\cU_q(\mathfrak{g})$-module $\cF(\Lambda_k).$
\end{thm}

\textbf{Note.}  If we replace the definition of $E_i$ and $F_i$ by
 \[
E_i=\sum_{\stackrel{(r,y) \in\bbN \times \bbZ}{k+y \equiv \pm i \bmod 2n}} \big (\prod_{\stackrel{(k_1,y_1)<(r,y)}{k_1+y_1 \equiv \pm i \bmod 2n}}T_{(k_1,y_1)}^- \big ) E_{(r,y)},
\]

\[
F_i=\sum_{\stackrel{(r,y) \in \times \bbN \times \bbZ}{k+y \equiv \pm i \bmod 2n}} \big (\prod_{\stackrel{(k_1,y_1)>(r,y)}{k_1+y_1 \equiv \pm i \bmod 2n}}T_{(k_1,y_1)}^+ \big ) F_{(r,y)},
\]

the pair $(L(\cF(\Lambda_k)),B(\cF(\Lambda_k)))$ is a (lower) crystal base of $\cF(\Lambda_k).$

\begin{defn} \label{sec3.2} 
Let $\my \in \cY(\Lambda_k)$ and $(k_1,y_1) >(k_2,y_2)>\cdots >(k_l,y_l)$ be the sites of the $i$-coloured corners of $\my$ and, define the $i$-signature of $\my$ to be the $l$-tuple $\sigma=(\sigma_1,\ldots,\sigma_l),$ where for $1 \leq m \leq l,$  
\[
\sigma_m:=\begin{cases}
 0 &\mbox{ if there is a concave corner in } \my \mbox{ at site } (k_m,y_m)\\ 
1 &\mbox{ if there is a convex corner in } \my \mbox{ at site } (k_m,y_m).
\end{cases}.
\]

  Define 
$J(\sigma)$\index{$J(\sigma)$} as follows:
let $J = \{ 1, \ldots, m \}$.
\begin{enumerate}
    \item[(i)]
    If there exists $r < s$ such that $(\sigma_{r}, \sigma_{s}) = 
    (0,1)$ and $r^{\prime} \not \in J$ for $r < r^{\prime} < s$, 
    replace $J$ by $J \backslash \{r, s\}$ and repeat this step;
    \item[(ii)]
    otherwise let $J(\sigma) = J$.
\end{enumerate}

If there exists an $i_r \in J(\sigma)$ with $\sigma_{i_r}=1,$
define $\etii(\my)$ to be the same as $\my$ with the $i$-coloured convex corner corresponding to the largest such $i_r $ removed; otherwise, define $\etii(\my)$ to be zero.

If there exists an $i_r \in J(\sigma)$ with $\sigma_{i_r}=0,$
define $\ftii(\my)$ to be the same as $\my$ with the $i$-coloured concave corner corresponding to the smallest such $i_r $ replaced by an $i$-coloured convex corner; otherwise, define $\ftii(\my)$ to be zero.

\end{defn}

\begin{thm}
The operators $\etii \mas \ftii$ defined above coincide with the Kashiwara's operators (see \cite{Ka91} for the definition of the Kashiwara's operators).
\end{thm}

\vfill
\newpage

\setlength{\unitlength}{0.20 cm}
\begin{center}
	\begin{picture}(80, 110)(0,0)
	\put(0,110){The crystal graph $ B(\Lambda_0)$ for $C_2^{(1)}.$}
	\put(49,100){$\phi$}
	\put(49.5,99) {\vector(0,-1){6}} 
	\put(50, 97){\tiny{0}}
	\put(46, 89){\begin{Young}0\cr\end{Young}}
	\put(49, 88){\vector(-1,-1){6}}
	\put(45, 86){\tiny{1}}
	\put(38, 78){\begin{Young}0&1\cr \end{Young}}
	\put(37, 77){\vector(-1,-1){7}}
	\put(47, 77){\vector(1,-1){7}}
	\put(33, 75){\tiny{1}}
	\put(50, 75){\tiny{2}}
	\put(22, 65){\begin{Young}0 & 1\cr\ 1 \cr \end{Young}}
	\put(51, 67){\begin{Young}0 & 1 & 2\cr \end{Young}}
	\put(25, 64){\vector(0,-1){11}}
	\put(24, 60){\tiny{0}}
	\put(26, 64){\vector(1, -1){11}}
	\put(32, 60){\tiny{2}}
	\put(59, 65){\vector(-1, -1){12}}
	\put(52, 60){\tiny{1}}
	\put(20, 47){\begin{Young}0 & 1\cr 1 &  0 \cr \end{Young}}
	\put(32, 47){\begin{Young}0 & 1 & 2 \cr 1 \cr \end{Young}}
	\put(42, 49.5){\begin{Young}0 & 1& 2&1 \cr \end{Young}}
	\put(25, 46){\vector(1,-1){10}}
	\put(26, 42){\tiny{2}}
	\put(38, 47){\vector(0, -1){11.5}}
	\put(36.5, 38){\tiny{0}}
	\put(41, 48.5){\vector(1, -1){13.5}}
	\put(52, 38){\tiny{2}}
	\put(44, 49){\vector(-1,-1){13.5}}
	\put(35, 41){\tiny{1}}
	\put(48, 47){\vector(0, -2){11.5}}
	\put(46.5, 38){\tiny{0}}
	\put(19, 30){\begin{Young}0 & 1 & 2 & 1\cr 1 \cr \end{Young}}
	\put(30.5, 30){\begin{Young}0 & 1 & 2 \cr 1& 0 \cr \end{Young}}
	\put(39, 32){\begin{Young}0 & 1& 2 & 1 & 0 \cr \end{Young}}
	\put(53, 28){\begin{Young}0 & 1& 2 \cr 1 \cr 2 \cr \end{Young}}
	\put(20, 20){$\vdots$}
	\put(35, 20){$\vdots$}
	\put(50, 20){$\vdots$}
	\put(65, 20){$\vdots$}

	\end{picture}
	

\end{center}

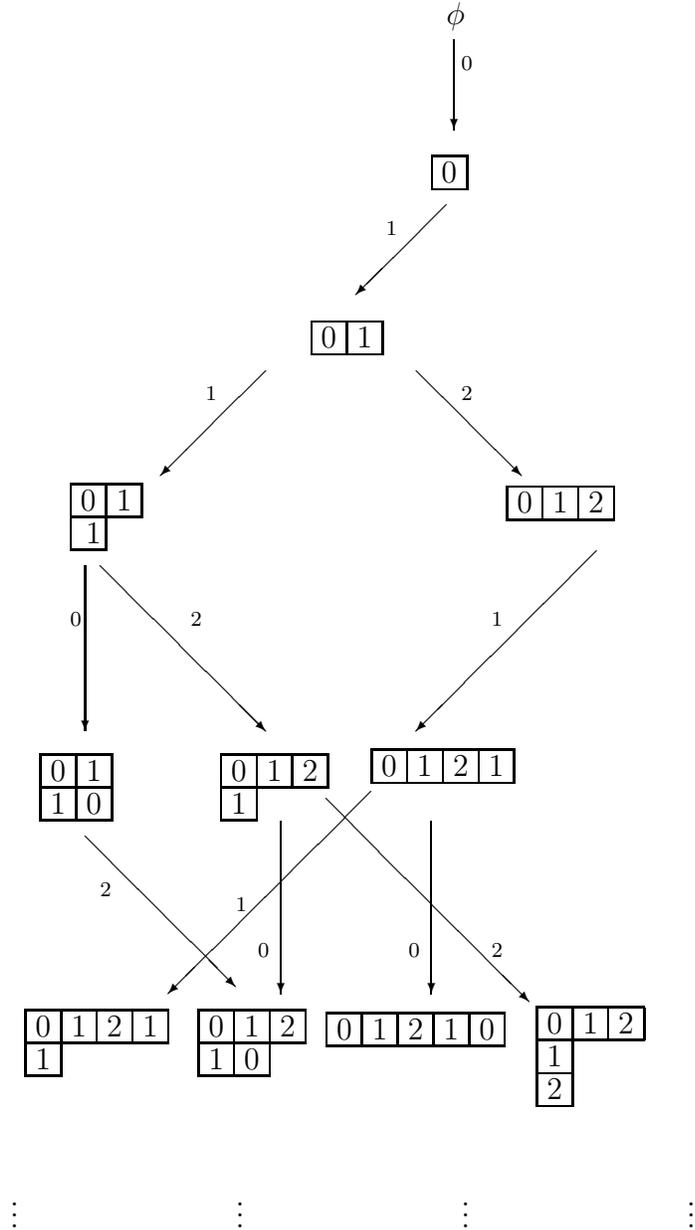
\setlength{\unitlength}{0.25 cm}
\begin{center}
	\begin{picture}(60,90)(0,0)
	\put(0,90){ The crystal graph $\ B(\Lambda_1)$ for $C_2^{(1)}.$}
	\put(30,80){$\phi$}
	\put(30.5,79) {\vector(0,-1){2}} 
	\put(31, 78){\tiny{1}}
	\put(28, 74){\begin{Young}1\cr\end{Young}}
	\put(27, 73){\vector(-1,-1){3}}
	\put(25, 72){\tiny{0}}
	\put(33, 73){\vector(1,-1){3}}
	\put(35, 72){\tiny{2}}
	\put(20, 65){\begin{Young}1\cr0\cr \end{Young}}
	\put(35, 65){\begin{Young}1 & 2\cr \end{Young}}
	\put(25, 63){\vector(1,-1){2}}
	\put(26, 62.5){\tiny{2}}
	\put(35, 63){\vector(-1,-1){2}}
	\put(33.5, 62.5){\tiny{0}}
	\put(22,63) {\vector(0,-1){2}} 
	\put(21,61.5){\tiny{1}}	
	\put(39,63) {\vector(0,-1){2}} 
	\put(40, 61.5){\tiny{1}}
	\put(19, 55){\begin{Young}1\cr 0\cr 1 \cr \end{Young}}
	\put(27, 56){\begin{Young}1 & 2\cr 0 \cr\end{Young}}
	\put(36, 58){\begin{Young}1 & 2 & 1\cr \end{Young}}
	\put(30, 54){\vector(0,-1){3}}
	\put(30.5, 52.5){\tiny{1}}
	\put(22, 54){\vector(1,-1){2}}
	\put(23, 53){\tiny{2}}
	\put(42, 57){\vector(-1,-1){5}}
	\put(39, 55){\tiny{0}}
	\put(20, 45){\begin{Young}1&2\cr 0\cr 1 \cr \end{Young}}
	\put(27, 47){\begin{Young}1&2&1\cr 0\cr \end{Young}}
	\put(34, 49){\begin{Young}1&2&1&0\cr \end{Young}}	
	\put(30, 46){\vector(0,-1){6}}
	\put(30.5, 42){\tiny{1}}
		\put(24, 44){\vector(0,-1){4}}
	\put(24.5, 41.5){\tiny{1}}
	\put(25, 46){\vector(4,-1){20}}
	\put(40.5, 42.5){\tiny{2}}
		\put(40, 47.5){\vector(-4,-1){25}}
	\put(18, 42.5){\tiny{0}}
		\put(38, 47.5){\vector(0,-1){7}}
	\put(38.5, 44.5){\tiny{1}}
	\put(10, 36){\begin{Young}1&2&1&0\cr 0\cr \end{Young}}
	\put(20, 34){\begin{Young}1&2&1\cr 0 \cr 1\cr \end{Young}}
	\put(27, 36){\begin{Young}1&2&1\cr 0&1\cr \end{Young}}	
	\put(34, 37){\begin{Young}1&2&1&0&1\cr \end{Young}}
	\put(45, 33){\begin{Young}1&2\cr 0\cr 1 \cr 2 \cr \end{Young}}	
	\put(15, 30){$\vdots$}
	\put(30, 30){$\vdots$}
	\put(45, 30){$\vdots$}
	
	\end{picture}
	

\end{center}

}

\end{document}